\documentclass[12pt,reqno]{amsart}
\usepackage{amsmath, amsfonts, amssymb, amsthm, amsbsy}

\theoremstyle{plain}
\newtheorem{theor}{Theorem}
\newtheorem{proposition}[theor]{Proposition}

\theoremstyle{definition}

\newtheorem{define}{Definition}
\newtheorem{example}{Example}

\theoremstyle{remark}
\newtheorem{rem}{Remark}

\newcommand{\BBR}{{\mathbb{R}}}

\DeclareMathOperator{\conv}{conv}

\begin{document}

\title{On eligibility by the de Borda voting rules}
\date{February 24, 2006, revised June 26, 2006}
\author{V.~Yu.~Kiselev}
\address{Department of Higher Mathematics,
Ivanovo State Power University,
Rabfakovsaya 34, Ivanovo 153003, Russia.}
\email{vkiselev@math.ispu.ru}


\begin{abstract}
We show that a necessary condition
for eligibility of a candidate by the set of de Borda's
voting rules in~[{H.~Moulin}
(1988) \textit{Axioms of cooperative decision making}]
is not sufficient
and we obtain some criterion for eligibility.
Let $r(a_i)$~be the score vector of a candidate~$a_i$,
$R$~be the set of all vectors $r(a_i)$,
and let $R'$~be the Pareto boundary of the convex hull~$\conv R$.
Then there is a scoring~$s$ such that
a candidate~$a$ wins with respect to the de Borda voting rule~$\beta_s$
if and only if~$r(a)\in R'$. 
 
\end{abstract}

\maketitle

\subsection*{Introduction}
Suppose the set~$A$ of candidates and the profile~$u$ of the voters'
preferences are fixed.
Let~$s$ be a system of scores and~$\beta_s$ be the
de Borda rule assigned to~$s$.
Further, let $\beta_s(u)\subset A$ be the set of winners
w.r.t.\ this rule.
A candidate~$a\in A$ is \emph{eligible}
w.r.t.\ the set of de Borda's rules~$\beta_s$
if there is a scoring~$s$ such that~$a\in\beta_s(u)$.
The book~\cite{Moulin} suggests the following necessary condition
for eligibility of a given candidate~$a$ w.r.t.\ the set of de Borda's rules:
a winner~$a$
has the score vector~$r(a)$ that belongs to the Pareto boundary of
the set~$R=\{r(a_j)$,\ $a_j\in A\}$;
see pages~\pageref{deB}\,--\,\pageref{P-boundary}
for rigorous definitions.

Unfortunately, the condition in this formulation 
is not sufficient. On page~\pageref{ShapEx} we give a counter\/-\/example
using the profile~$u$ defined in Eq.~\eqref{Shapovalov}.

In Theorems~\ref{RPO} and~\ref{SndTh} below
we prove some version of the eligibility criterion and so we give a complete
solution for a problem from~\cite[p.~249]{Moulin}.
The difference between the condition in our criterion and the necessary
condition from~\cite{Moulin}
is that the Pareto boundary of the \emph{convex hull}~$\conv R$
must be used instead of the Pareto boundary of the set~$R$ itself,
here~$R$ is the set of all score vectors.

\section{Definitions}
\subsection{Profiles}
Let $P_1$,\ $\dots$,\ $P_n$ be \emph{electors} (or
\emph{voters}) and $A=\{a_1$;\ $\dots$;\ $a_p\}$ be the set of
\emph{candidates} in some elections.
Suppose that every elector $P_i$ has an opinion about each
candidate such that the candidates
are arranged by the strict order~${}>_i{}$:
the first candidate in this rearrangement is the most favourable
for~$P_i$, etc. This strict linear order ${}>_i{}$ on~$A$ is
called the \emph{preference} of the elector~$P_i$ and is
denoted by~$u_i$. The order~$u_i$ is given by the sequence
\[
a_{j_1}\mathrel{>_i}a_{j_2}\mathrel{>_i}\dots\mathrel{>_i}a_{j_p},
\]
where $J=(j_1;j_2;\dots;j_p)$ is a rearrangement of $(1;2;\dots;p)$;
generally, $J$~depends on the elector~$P_i$.
%

In what follows, 
we write down the elements of the preferences~$u_i$ in columns
and thus we compose the matrix  
\[
u=\bigl|u_1\ u_2\ \dots\ u_n\bigr|.
\]
This matrix is the \emph{profile} of preferences of the voters.
Let us denote by~$L(A)$ the set of all strict linear orders on~$A$.
Then the profiles~$u$ are elements of the \emph{space of
preferences}~${\mathcal{L}}=(L(A))^n$.

Suppose that several voters have coinciding preferences~$u_i$,
and assume that the order of electors
is not important for the voting rules.
Then we write down the coinciding columns only once and
indicate the respective number of
electors in the upper line of the matrix~$u$, e.~g.,
\begin{equation}\label{Shapovalov}
u=\begin{vmatrix}
{}^2&{}^6&{}^7&{}^1\\
a&b&c&c\\
b&a&a&b\\
c&c&b&a
\end{vmatrix}.
\end{equation}
Notation~\eqref{Shapovalov} means that
$p=3$,\ $n=16$, and $A=\{a;b;c\}$ is the set of candidates; two
electors have the preference ${a>_ib>_ic}$ (${i=1,2}$), six electors
have the preference ${b>_ia>_ic}$ (${i=3,\dots,8}$), seven electors
have the preference ${c>_ia>_ib}$ (${i=9,\dots,15}$), and the $16$th
elector has the preference ${c>_{16}b>_{16}a}$.

We shall use the profile defined by Eq.~\eqref{Shapovalov}
in Example~\ref{ShapEx} on page~\pageref{ShapEx}
as a counter\/-\/example to the erroneous criterion in~\cite{Moulin}.

\begin{define}\label{def_elections}
Let $n$ and $p$ be positive integers.
A \emph{voting rule}~$\varphi$ is a point\/-\/to\/-\/set
map\footnote{A \emph{point\/-\/to\/-\/set map} $f\colon D\to E$
maps each $x\in D$ to some subset $f(x)\subset E$.}
${\varphi\colon{\mathcal{L}}\to A}$ such that $\varphi(u)\ne\varnothing$ for
each profile $u\in{\mathcal{L}}$ of dimension $n\times p$.
The candidates $a\in\varphi(u)$
are called \emph{winners} with respect to the rule~$\varphi$ and
the profile~$u$.
\end{define}

\subsection{Generalized de Borda's voting rules}
A generalized de Borda's voting rule is defined by the
following conditions.

\begin{define}\label{scores}
Let $s=(s_0;s_1;\dots;s_{p-1})$ be a vector such that
\begin{equation}\label{nestrogo}
s_0\le s_1\le \dots\le s_k\le \dots\le s_{p-1}\text{ and }
s_0<s_{p-1}.
\end{equation}
The vector~$s$ is a \emph{system of scores}, or a \emph{scoring}.
Then, each candidate $a\in A$ obtains $s_0$
points for the last place in the preference of a voter,
$s_1$~points for the last but one place, $s_2$~points for the
third place from the bottom, etc.; clearly,
$s_{p-1}$ points are obtained for the top position.
The points accumulated by a candidate~$a$ from all the voters are summed,
and this sum~$B_s(a)$ is called the \emph{de Borda estimate} of the
candidate~$a$ w.r.t.\ the scoring~$s$, or, briefly, the
\emph{$s$-estimate} for~$a$.
The \emph{de Borda voting rule}~$\beta_s$\label{deB}
with the scoring~$s$ says that
the candidates $a\in A$ who have the maximal sum $B_s(a)$ of points
are the \emph{winners} of these elections.
Note that the winner may be not unique
if two or more candidates $a,b,\dots$ obtain equal
(maximal) estimates $B_s(a)=B_s(b)=\dots$\;.
We denote by~$\beta_s(u)$ the set
of all winners with respect to the profile~$u$.
\end{define}

\begin{example}
The \emph{standard de Borda voting rule} is based on the system of scores
$s=(0;1;2;\dots;p-1)$: a candidate obtains $0$ points
for the last place, $1$ point
for the last but one place, etc.; a candidate obtains $p-1$ points for
the first place. The rule suggests that the winner has the maximal sum
of points from all voters.

The \emph{plurality voting rule} 
is also a generalized de Borda's rule. The rule corresponds to
the scoring~$s$ with $s_0=0$,\ $s_1=0$,\ $\dots$,\
$s_{p-2}=0$,\ $s_{p-1}=1$. Therefore
the sum~$B_s(a)$ consists of~$m$ units if $a$~is the preferable
candidate for $m$~voters;
hence the winner by the rule $\beta_s$ receives the plurality
in the preferences of the voters.

The \emph{antiplurality voting rule} 
is another generalized de
Borda's rule. It corresponds to the scoring~$s$ such that
$s_0=-1$, $s_1=0$, $s_2=0$, $\dots$, $s_{p-1}=0$.  The candidate who
wins by the rule $\beta_s$ has the least number of the
last places in the preferences of the electors.
\end{example}

\begin{rem}\label{s0}
Let $\beta_s$ be a de Borda's voting rule with the scoring
$s=(s_0;s_1;\dots;s_{p-1})$.
If a constant $d$ is added to each score~$s_k$, $k=0,\dots,p-1$,
then we get the new rule $\beta_{s'}$ with $s'$
such that ${s'_k=s_k+d}$. Hence all the estimates~$B_s(a)$
change to ${B_{s'}(a)=B_s(a)+nd}$, that is, they acquire the
addend $nd$, which does not depend on a candidate~$a$. Obviously,
the sets of winners by the two rules~$\beta_s$ and~$\beta_{s'}$
coincide: ${\beta_{s'}(u)=\beta_s(u)}$ for any profile~$u$. Consequently,
the scorings~$s$ and~$s'$ define \emph{the same} voting rule:
${\beta_{s'}=\beta_s}$. If we set ${d=-s_0}$, then we have ${s_0'=0}$.
Therefore we can 
set $s_0=0$ for the de Borda rule assigned to any
scoring~$s$ whenever that is convenient.
\end{rem}

\begin{define}
Let $F\subset\BBR^m$ be a closed set. Let $f'\in F$ be a point such that
there are no other points $f\in F$ which satisfy the following two
conditions
\begin{itemize}
\item
$f_j\ge f'_j$ for all $j=1,\dots,m$, and
\item
$f_{j_0}>f'_{j_0}$ for some $j_0$.
\end{itemize}
The \emph{Pareto boundary}\label{P-boundary}
${F'\subset F}$ is the subset constituted by all the points~$f'$.

A point $f''\in F$ belongs to the
\emph{weak Pareto boundary} $F''$ of the set $F$ 
if there are no points $f\in F$ such that all
coordinates of~$f$ are greater than
the respective coordinates of~$f''$.
\end{define}

Finally, let us recall a helpful statement from the convex analysis.

\begin{proposition}\label{prop1}
Let~$F$ be a bounded convex subset of $\BBR^m$ and $F'$~be the
\textup{(}weak\textup{)} Pareto boundary of~$F$.
Then for any point $f^0\in F'$
there is a linear function $\ell(f)=d_1f_1+\dots+d_mf_m$ with positive
\textup{(}resp.\textup{,} nonnegative\textup{)} coefficients $d_j$
such that $\ell$ achieves its
maximum on~$F$ at~$f^0$. The converse is also true\textup{:} any point
$f^{00}\in F\setminus F'$ is not a point of maximum on~$F$ for any
linear function with positive
\textup{(}resp.\textup{,} nonnegative\textup{)} coefficients.
\end{proposition}

\section{Eligible and uneligible candidates}
Let a profile~$u$ of the preferences be fixed. Suppose that the winners are
determined by a de Borda's rule~$\beta_s$ from some fixed set $\boldsymbol{\beta}$
of these rules,
but we do not know which particular scoring~$s$ will be used. Now we want
to predict which candidates $a\in A$ can win by a rule~$\beta_s\in\boldsymbol{\beta}$
assigned to some scoring~$s$ (these candidates are
\emph{eligible} w.r.t.\ this set of rules $\boldsymbol{\beta}$)
and which candidates definitely can not win for any~$\beta_s\in\boldsymbol{\beta}$ (we call
them \emph{uneligible} w.r.t.\ this set of rules).

We claim that the set of eligible candidates can be found in two important
cases (see Theorems~\ref{RPO} and~\ref{SndTh}) without knowing beforehand the
scoring~$s$.

Suppose a profile~$u$ is fixed and~$a$ is a candidate.
The following procedure allows to obtain
the \emph{score vector}~$r(a)$ for~$a$ using the profile~$u$.
By construction, the score vector~$r(a)$ has $p-1$ components.
The first component~$r_1(a)$ equals
the number of electors who regard~$a$ as the best. The second
component~$r_2(a)$ is equal to
the number of voters for whom $a$~has either the first
or the second place in their preferences; the component $r_3(a)$ is
equal to the number of electors having~$a$ on some of three first
places, etc.                      
Thus we obtain the score vector
${r(a)=(r_1(a);r_2(a);\dots;r_{p-1}(a))}$ for every candidate $a\in A$.
Obviously, the components of the vector~$r(a)$ are non\/-\/decreasing:
$$
r_1(a)\le r_2(a)\le\dots\le r_{p-1}(a).
$$
Also, we note that the score vector~$r(a)$ is found without knowing the
scoring~$s$.

\label{Rvect}Consider the set $R\subset\BBR^{p-1}$ that consists of
$p$ score vectors~$r(a)$ for all candidates $a\in A$. Next, find
the Pareto boundary~$R'$ of the convex hull $\conv R$ of~$R$.

Now we formulate the main result of this article. Theorem~\ref{RPO}
contains the criterion
of eligibility of a given candidate by de Borda's rules with strict scorings.

\begin{theor}\label{RPO}
Suppose a profile~$u$ is fixed. Consider the set $\boldsymbol{\beta}$ of all de
Borda's rules $\beta_s$ with the strict scorings~$s$ such that
\begin{equation}\label{strogo}
s_0<s_1<\dots<s_{p-1}.
\end{equation}
Let $a\in A$ be a candidate and $r(a)$ be the score vector for~$a$.
Then the candidate~$a$ is eligible w.r.t.\ the set of rules $\boldsymbol{\beta}$
if and only if $r(a)\in R'$, where
$R'$~is the Pareto boundary of the convex hull $\conv R$ of the set
$R=\{r(a_j)$,\ $a_j\in A$,\ $j=1\ldots p\}$.
\end{theor}

\begin{proof}
According to Remark~\ref{s0}, we suppose that $s_0=0$.
Let us show that if $a\in\beta_s(u)$ with some scoring~$s$
such that inequalities~\eqref{strogo} hold, then $r(a)\in R'$.
Let the differences between the scores be
${d_1=s_1-s_0=s_1}$, ${d_2=s_2-s_1}$,
$d_3=s_3-s_2$, $\dots$, $d_{p-1}=s_{p-1}-s_{p-2}$.
We note that~\eqref{strogo} implies ${d_{p-k}>0}$ for all
${k=1,2,\dots,p-1}$.

Now we analyze the summands that contribute to the de Borda estimate~$B_s(a)$.
The candidate~$a$ has the first place
in preferences of $r_1(a)$~voters and hence obtains ${r_1(a)\cdot s_{p-1}}$
points. Next,
${r_2(a)-r_1(a)}$ voters put the candidate~$a$ on the second place,
therefore $a$~obtains ${(r_2(a)-r_1(a))\cdot s_{p-2}}$ more points.
Similarly, for the $k$th place in
the preferences of ${(r_k(a)-r_{k-1}(a))}$ voters
the candidate~$a$ obtains
${(r_k(a)-r_{k-1}(a))\cdot s_{p-k}}$ points in the sum~$B_s(a)$, here
${k=1,2,\dots,p-1}$. Recall that the last place results in no points
since~${s_0=0}$.

By construction, the estimate~$B_s(a)$ is the sum of points obtained by~$a$:
\begin{multline}\label{Bdiff}
B_s(a)=
r_1(a)s_{p-1}+(r_2(a)-r_1(a))s_{p-2}+(r_3(a)-r_2(a))s_{p-3}+\dots+\\
+(r_k(a)-r_{k-1}(a))s_{p-k}+\dots+(r_{p-1}(a)-r_{p-2}(a))s_1=\\
=r_1(a)(s_{p-1}-s_{p-2})+r_2(a)(s_{p-2}-s_{p-3})+\dots+r_{p-1}s_1=\\
=r_1(a)d_{p-1}+r_2(a)d_{p-2}+\dots+r_{p-1}d_1
=\sum_{k=1}^{p-1}d_{p-k}r_k(a).
\end{multline}
Thus we obtain the formula for the de Borda estimate with the scoring~$s$
(see~\cite[Ch.~9]{Moulin}). For $a\in\beta_s(u)$,
the sum in~\eqref{Bdiff} is maximal and the vector~$r(a)$ is a
point of maximum on~$R$ for the linear function
${\ell(r)=\sum\limits_{k=1}^{p-1}d_{p-k}r_k}$
with positive coefficients~$d_{p-k}$.   
Moreover, we note that the maximum of a linear function on the convex hull of
a set equals the maximum of the function
on the set itself. Hence, according to
Proposition~\ref{prop1}, the maximum of the function~$\ell(r)$
is achieved only at a point of the set~$R'$, that is,
$r(a)\in R'$ if $a\in\beta_s(u)$.

Conversely, suppose $r(a)\in R'\cap R$. By Proposition \ref{prop1},
we can find a system of positive coefficients~$d_{p-k}$
such that the linear function~$\ell(r)$ with these coefficients has a maximum
on~$R'$ (and, consequently, on~$R$) at the point~$r(a)$.
Using these coefficients, we construct the scoring
$s=(s_0;s_1;\dots;s_{p-1})$ by setting $s_0=0$, $s_1=d_1$, $s_2=d_1+d_2$,
$s_3=d_1+d_2+d_3$, $\dots$. Then from Eq.~\eqref{Bdiff} it follows
that the sum~$B_s(a)$ achieves its maximum at~$a\in A$
and therefore~$a\in\beta_s(u)$.
\end{proof}

Obviously, we can modify slightly the assumptions
of Theorem~\ref{RPO} and consider a wider set
of de Borda's rules, namely, the set $\boldsymbol{\beta}'$ of rules~$\beta_s$ with
the scorings~$s$ that satisfy condition~\eqref{nestrogo}
but may not satisfy~\eqref{strogo}.
Then, using Proposition~\ref{prop1} again, we obtain

\begin{theor}\label{SndTh}
Fix a profile~$u$ and consider the set $\boldsymbol{\beta}'$
of all de Borda's rules~$\beta_s$
with the scorings~$s$ that satisfy condition~\eqref{nestrogo}.
Find the score vectors~$r(x)$ for all candidates $x\in A$
and consider the set~$R$ that consists of all vectors~$r(x)$.
Let $R''$~be the weak Pareto boundary of the convex hull of~$R$.
Then any candidate~$a$ such that $r(a)\in R''$ is eligible
w.r.t.\ this set of voting rules, and only these candidates are eligible.
\end{theor}

\begin{rem}
Of course,
the conclusions of the two theorems are not true if the set
$R\cap R'$ (or $R\cap R''$)
is replaced by the Pareto boundary (or the weak Pareto boundary)
of the set~$R$ itself.
(We recall that the set $R\cap R'$ (or $R\cap R''$)
consist of the vectors~$r(a)$
that belong to the Pareto boundary (or to the weak Pareto boundary)
of~$\conv R$.)
Indeed, the sets $R\cap R'$ and $R\cap R''$ can be more narrow
than the (weak) Pareto boundary of~$R$.\footnote{%
   Clearly,
   since $R\cap R'$ is a part of the Pareto boundary of~$R$ (and
   $R\cap R''$ is a part of the weak Pareto boundary of $R$),
   the condition of~\cite{Moulin} is necessary.
   Namely, any winner~$a$ by some rule~$\beta_s$ has
   the (weak) Pareto optimal (i.~e., belonging
   to the (weak) Pareto boundary of~$R$) score vector~$r(a)$.
   But this condition is not sufficient:
   not every candidate~$a$ having the (weak) Pareto optimal score vector~$r(a)$
   can win by some voting rule~$\beta_s$.}
This situation is illustrated by the following example.
\end{rem}

\begin{example}[A.~V.~Shapovalov, private communication]\label{ShapEx}   
Let $p=3$ and $n=16$. Consider profile~\eqref{Shapovalov};
then we have ${r(a)=(2;15)}$, ${r(b)=(6;9)}$,
${r(c)=(8;8)}$.

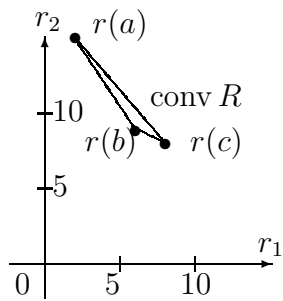
\begin{figure}[h]
\begin{center}
\unitlength 1mm
\linethickness{0.4pt}
\begin{picture}(35.00,37.33)
\put(5.00,35.00){\vector(0,1){0.2}}
\put(5.00,0.33){\line(0,1){34.67}}
\multiput(15,4)(10,0){2}{\line(0,1){2}}
\put(35.00,5.00){\vector(1,0){0.2}}
\put(0.33,5.00){\line(1,0){34.67}}
\multiput(4,15)(0,10){2}{\line(1,0){2}}
\multiput(9.00,35.00)(0.12,-0.18){67}{\line(0,-1){0.18}}
\multiput(17.00,23.00)(0.24,-0.12){17}{\line(1,0){0.24}}
\multiput(21.00,21.00)(-0.12,0.14){100}{\line(0,1){0.14}}
\put(9.00,35.00){\circle*{1.33}}
\put(17.00,22.67){\circle*{1.33}}
\put(21.00,21.00){\circle*{1.33}}
\put(1.00,1.00){\makebox(0,0)[lb]{$0$}}
\put(13,1){\makebox(0,0)[lb]{$5$}}
\put(23,1){\makebox(0,0)[lb]{$10$}}
\put(6,14){\makebox(0,0)[lb]{$5$}}
\put(6,24){\makebox(0,0)[lb]{$10$}}
\put(10.33,19.00){\makebox(0,0)[lb]{$r(b)$}}
\put(19,26){\makebox(0,0)[lb]{$\conv R$}}
\put(24.33,19.00){\makebox(0,0)[lb]{$r(c)$}}
\put(11.33,34.67){\makebox(0,0)[lb]{$r(a)$}}
\put(3.67,36){\makebox(0,0)[lb]{$r_2$}}
\put(33.33,6){\makebox(0,0)[lb]{$r_1$}}
\end{picture}
\end{center}
\caption{Three vectors $r(a)$, $r(b)$, $r(c)$ and their convex hull
$\conv R$}\label{fig1}
\end{figure}

The Pareto boundary of the set ${R=\{r(a);r(b);r(c)\}}$
contains all the three vectors $r(a)$, $r(b)$, and~$r(c)$,
but only~$r(a)$ and~$r(c)$ belong to the sets~$R'$ and $R''$
(see Fig.~\ref{fig1}).

Hence only~$a$ and~$c$ can be the winners
by some of the generalized de Borda rules if we have
profile~\eqref{Shapovalov}. The vector~$r(b)$
is Pareto optimal, but none of the linear functions
${\ell(r_1;r_2)=d_2r_1+d_1r_2}$
with positive or nonnegative
coefficients~$d_1$ and~$d_2$ achieves its maximum on~$R$
at the point~$r(b)$. Therefore the candidate~$b$ is uneligible neither w.r.t.\
the set $\boldsymbol{\beta}$ nor w.r.t.\ the set $\boldsymbol{\beta}'$:
there is no rule~$\beta_s$ such that~$b$ is a winner.
\end{example}

\medskip

The author thanks Prof.\ A.~V.~Shapovalov (KTH, Stockholm)
for illuminating discussions and Prof.\ H.~Moulin (Rice University,
Houston, Texas) for favourable attention to the article and the
approval of its results.

\end{document}